\documentclass[letterpaper, 10 pt, conference]{ieeeconf}  

\IEEEoverridecommandlockouts
\usepackage[T1]{fontenc}
\usepackage{microtype}

\usepackage{xcolor}

\usepackage{amsmath}                                      
\usepackage{amssymb}                                      
\usepackage{bm}                                           
\usepackage{commath}                                      
\usepackage{mathtools}                                    
\usepackage{textcomp}
\usepackage{siunitx}                                      

\usepackage{booktabs}                                     
\usepackage{cite}                                         

\makeatletter
\let\NAT@parse\undefined
\makeatother
\usepackage{hyperref}

\usepackage{tikz}
\usepackage{pgfplots}
\usepackage{pgfplotstable}

\pgfplotsset{
    compat=1.16,
    cycle list name=linestyles*,
    every axis plot/.append style={
        thick
    },
    log x ticks with fixed point/.style={
        xticklabel={
            \pgfkeys{/pgf/fpu=true}
            \pgfmathparse{exp(\tick)}%
            \pgfmathprintnumber[fixed relative, precision=3]{\pgfmathresult}
            \pgfkeys{/pgf/fpu=false}
        }
    },
    log y ticks with fixed point/.style={
        yticklabel={
            \pgfkeys{/pgf/fpu=true}
            \pgfmathparse{exp(\tick)}%
            \pgfmathprintnumber[fixed relative, precision=3]{\pgfmathresult}
            \pgfkeys{/pgf/fpu=false}
        }
    }
}

\allowdisplaybreaks                                        

\newtheorem{assumption}{Assumption}
\newtheorem{theorem}{Theorem}
\newtheorem{lemma}[theorem]{Lemma}

\newcommand{\name}[1]{\texttt{#1}}

\newcommand{\bbma} {\begin{bmatrix} }
\newcommand{\ebma} {\end{bmatrix}}
\newcommand{\real} {\mathbb{R}}

\DeclareMathOperator{\diag}{diag}

\newcommand{\correction}[1]{\textcolor{blue}{#1}}

\usepackage{fancyhdr}

\fancypagestyle{firstpage}{%
    \fancyfoot[C]{\small
        © 2021 IEEE.
        Personal use of this material is permitted.
        Permission from IEEE must be obtained for all other uses, in any current or future media, including reprinting/republishing this material for advertising or promotional purposes, creating new collective works, for resale or redistribution to servers or lists, or reuse of any copyrighted component of this work in other works.%
    }%
}

\title{\LARGE \bf
Convergence Properties of Fast quasi-LPV Model Predictive Control
}

\author{Christian Hespe$^{1}$ and Herbert Werner$^{1}$
\thanks{$^{1}$Hamburg University of Technology, Institute of Control Systems, Eissendorfer Str. 40, 21073 Hamburg, Germany. \texttt{\{christian.hespe, h.werner\}@tuhh.de}}%
}

\begin{document}

    \maketitle
    \thispagestyle{firstpage}

    \begin{abstract}
        In this paper, we study the convergence properties of an iterative algorithm for fast nonlinear model predictive control (MPC) of quasi-linear parameter-varying systems without inequality constraints.
Compared to previous works considering this algorithm, we contribute conditions under which the iterations are guaranteed to converge.
Furthermore, we show that the algorithm converges to suboptimal solutions and propose an optimality-preserving variant with moderately increased computational complexity.
Finally, we compare both variants in terms of quality of solution and computational performance with a state-of-the-art solver for nonlinear MPC in two simulation benchmarks.

    \end{abstract}

    \pagestyle{empty}
    \section{Introduction}\label{sec:introduction}
Model predictive control (MPC) is an advanced control strategy that calculates its control trajectory by solving a sequence of optimization problems.
Because of its ability to handle constraints in an optimal fashion as well as straight-forward treatment of multi-input multi-output and nonlinear plants, MPC can be applied to broad classes of control problems.
However, the application of MPC is often hindered by its extensive need for computational resources caused by solving potentially nonconvex and nonlinear optimization problems online.
For this reason, MPC was first applied in the process industry, where plants are often characterized by slow dynamics and long sampling times \cite{Rawlings2017}.

In recent years, the application of MPC has become feasible for plants with increasingly fast dynamics.
This is due in part to exponential growth in computing power, but is likewise enabled by advances in algorithms for solving nonlinear programs (NLPs) and development of fast software for nonlinear MPC.
Examples for such software include \name{acados} \cite{Verschueren2019}, which is based on sequential quadratic programming (SQP) and real-time iterations as proposed in \cite{Diehl2005}, as well as \name{GRAMPC} \cite{Englert2019} and \name{OpEn} \cite{Sopasakis2020}, which rely on first-order augmented Lagrangian methods.
Further improvements of the computational performance can be achieved by approximating the solution of the optimal control problem (OCP), as it is done among others in zero-order MPC \cite{Zanelli2019a}.

In this paper, we will consider an algorithm for fast nonlinear MPC of quasi-linear parameter-varying (quasi-LPV) systems that was originally proposed in \cite{Cisneros2016}.
By fixing the parameter trajectory of the quasi-LPV system during optimization and thus temporarily disregarding the dependency of the parameter on states and inputs, the system is converted into a linear \emph{time}-varying model, for which MPC can be formulated as a quadratic program (QP).
This procedure is then iterated by calculating a new parameter trajectory from the solution and repeating the process.
The algorithm has been successfully applied to a variety of problems, c.f. \cite{Alcala2019, Cisneros2019, Gango2019, Hoekstra2019, MarquezRuiz2019, Mulders2019}.

Compared to standard nonlinear MPC, one advantage of the quasi-LPV framework is that it facilitates stability analysis.
Suitable terminal ingredients guaranteeing closed-loop stability for LPV models can be computed offline using linear matrix inequalities \cite{Cisneros2020}.
Additionally, due to the way the parameters are handled, the algorithm does not suffer from conservatism for models with a large number of parameters, in contrast to the approach of designing controllers robust with respect to the parameter trajectories.
This is utilized in \cite{Cisneros2020a}, where velocity-based linearization is applied to obtain LPV models with advantageous properties but many parameters.

The main contribution of this paper is an analysis of the algorithm's convergence properties, which was lacking from previous publications despite the heuristic evidence through successful applications.
For problems without inequality constraints, we provide conditions under which the quasi-LPV MPC (qLMPC) algorithm is guaranteed to converge and show that the iterations' fixpoints are suboptimal solutions to the original NLP.
An extension to the general case with both equality and inequality constraints is subject to future research.
Furthermore, it is shown that qLMPC is equivalent to SQP and thus calculates optimal solutions if a modified LPV model with increased computational complexity is used.
In order to assess the degree of suboptimality of the original variant and the impact of the complexity increase, we benchmark both qLMPC variants against \name{acados} and an optimal reference solution in two simulation scenarios.

The paper continues with a short review of the qLMPC algorithm in Section~\ref{sec:problem}.
Section~\ref{sec:analysis} contains the main result, the convergence analysis of the qLMPC algorithm.
Equivalence to SQP is shown in Section~\ref{sec:optimality}, and Section~\ref{sec:benchmarking} demonstrates the potential of the qLMPC algorithm in simulation benchmarks.
Finally, we present concluding remarks in Section~\ref{sec:conclusions}.

\section{Problem Setup}\label{sec:problem}
\subsection{Notation}\label{sec:problem_notation}
We denote the Kronecker product of matrices \(M_1\) and \(M_2\) with \(M_1 \otimes M_2\) and \(M \succ(\succeq)\ 0\) stands for the matrix \(M\) being symmetric positive (semi)definite.
Block diagonal matrices are written as \(M = \diag(M_1, \ldots, M_n)\), with elements \(M_1\) to \(M_n\) on the diagonal.
\(\|z\|_M \coloneqq \sqrt{z^T M z}\) with \(M \succ(\succeq)\ 0\) denotes the weighted vector (semi)norm and the notation \((z_1, \ldots, z_n) \coloneqq \left[ z_1^T, \ldots, z_n^T \right]^T\) is used for vertical concatenation.

\subsection{Quasi-LPV System Description}\label{sec:problem_system}
Consider the discrete-time LPV system described by
\begin{equation}\label{eq:lpv_system}
    x_{k+1} = A(\rho_k) x_k + B(\rho_k) u_k,
\end{equation}
where \(x_k \in \real^{n_x}\) is the state, \(u_k \in \real^{n_u}\) the input, \(\rho_k \in \mathbb{R}^{n_\rho}\) is a time-varying parameter and \(A(\rho): \mathbb{R}^{n_\rho} \to \real^{n_x \times n_x}\) and \(B(\rho): \mathbb{R}^{n_\rho} \to \real^{n_x \times n_u}\) are the parameter-varying model matrices of the system.
The initial condition is given by \(x_0 = \hat{x}_0\).
For quasi-LPV systems, \(\rho_k\) is a function of the state and input defined as
\begin{equation}\label{eq:lpv_rho_fun}
   \rho_k \coloneqq \rho(x_k, u_k): \real^{n_x} \times \real^{n_u} \to \mathbb{R}^{n_\rho}.
\end{equation}

We assume, without loss of generality, that the origin is an unforced equilibrium of the system, i.e. \(\bar{x} = A(\bar{\rho}) \bar{x} + B(\bar{\rho}) \bar{u}\) with \(\bar{x} = 0\), \(\bar{u} = 0\) and \(\bar{\rho} = \rho(\bar{x}, \bar{u})\).
Our goal is to solve the regulator problem, that is, to bring the system into this equilibrium.

\subsection{Quasi-LPV Model Predictive Control}\label{sec:problem_mpc}
We will control the plant defined by \eqref{eq:lpv_system} and \eqref{eq:lpv_rho_fun} using MPC.
Thus, the control input is calculated by solving a finite horizon OCP at each sampling instance, which for system~\eqref{eq:lpv_system} is given by
\begin{subequations}\label{eq:opt_lpv_sum}\begin{alignat}{2}
    &\!\min_{\bm{x}, \bm{u}} & \quad \|x_N&\|_P^2 + \sum_{k=0}^{N-1} \|x_k\|_Q^2 + \|u_k\|_R^2 \\
    &\text{s.t.} & x_0     &= \hat{x}_0 \label{eq:opt_lpv_sum_initial_value} \\
    &            & x_{k+1} &= A(\rho_k) x_k + B(\rho_k) u_k \quad \forall k \in [0, N-1] \label{eq:opt_lpv_sum_constraints_dynamics}
\end{alignat}\end{subequations}
where the state and input vectors of each step are collected in \(\bm{x} \coloneqq (x_0, x_1, \ldots, x_N)\) and \(\bm{u} \coloneqq (u_0, u_1, \ldots, u_{N-1})\) and \(Q, R, P \succeq 0\) are weighting matrices for the states, inputs and terminal states, respectively.
Because of the dependency of \(\rho_k\) on \(x_k\) and \(u_k\) this is not a quadratic but a nonlinear program.
Out of the solution vector~\(\bm{u}\), only \(u_0\) is applied to the plant, leading to receding horizon control.

To facilitate the analysis of the OCP, we reformulate the problem in a more general form.
By combining the states and inputs to \(\bm{y} \coloneqq (\bm{x}, \bm{u})\) and collecting the parameter vectors in \(\bm{\rho}(\bm{y}) \coloneqq (\rho_0, \rho_1, \ldots, \rho_{N-1})\), optimization problem~\eqref{eq:opt_lpv_sum} can be rewritten as
\begin{subequations}\label{eq:opt_cons}\begin{alignat}{2}
    &\!\min_{\bm{y}} & \qquad & \|\bm{y}\|_{\mathcal{Q}}^2\\
    &\text{s.t.}     &        & \mathcal{G}\left(\bm{\rho}(\bm{y})\right) \bm{y} + C \hat{x}_0 = 0 \label{eq:opt_cons_constraints_dynamics}
\end{alignat}\end{subequations}
with the new weighting matrix \(\mathcal{Q} \coloneqq \diag(I_N \otimes Q, P, I_N \otimes R)\), \(C \coloneqq \bbma -I_{n_x} & 0 & \cdots & 0\ebma^T\) and the matrix valued function \(\mathcal{G}(\bm{\rho})\).
In the following, the dependency of \(\bm{\rho}(\bm{y})\) on \(\bm{y}\) is often dropped for brevity of notation.
Note that if \eqref{eq:opt_lpv_sum} is reformulated as \eqref{eq:opt_cons}, then \(\mathcal{G}(\bm{\rho})\) is sparse and has full row rank for all values of \(\bm{\rho}\).
We will make the following assumptions during the convergence analysis in the next section.

\begin{assumption}\label{ass:model_differentiability}
    The functions \(\mathcal{G}(\bm{\rho})\) and \(\bm{\rho}(\bm{y})\) are twice continuously differentiable.
\end{assumption}

\subsection{Iterative Refinement}\label{sec:problem_iteration}
Instead of directly solving the NLP \eqref{eq:opt_cons}, it was first proposed in \cite{Cisneros2016} to solve the OCP in an iterative way akin to SQP.
In fact, as we will show in Section~\ref{sec:optimality}, a variant of this scheme is equivalent to SQP with a Gauss-Newton approximation.

The idea is to use a fixed parameter trajectory \(\bm{\rho}[l]\) while solving \eqref{eq:opt_cons} instead of letting \(\bm{\rho}(\bm{y})\) depend on \(\bm{y}\).
Here, \(l\) does not index through time like \(k\) but through iterations applied at a single time step.
In this way, we obtain a \emph{linear time}-varying plant model instead of the \emph{quasi-linear parameter}-varying one.
At every step, \(\bm{\rho}[l]\) is determined from \(\bm{y}[l]\), the solution of the previous iteration.
This leads to a sequence of QPs given by
\begin{subequations}\label{eq:opt_iter}\begin{alignat}{2}
    &\!\min_{\bm{y}} & \qquad & \|\bm{y}\|_{\mathcal{Q}}^2\\
    &\text{s.t.}     &        & \mathcal{G}(\bm{\rho}[l]) \bm{y} + C \hat{x}_0 = 0, \label{eq:opt_iter_constraints_dynamics}
\end{alignat}\end{subequations}
which can be solved efficiently using suitable algorithms.
For a detailed description of the approach see \cite{Cisneros2016, Cisneros2020}.

\section{Convergence Analysis}\label{sec:analysis}
The iterative procedure for solving nonlinear OCPs that was outlined in Section~\ref{sec:problem_iteration} is known to work well in practice.
However, no formal analysis of its convergence properties has been provided so far.
To obtain conditions for its convergence, we will proceed similar to the approach taken in \cite{Bock2007} by interpreting the iterative procedure as inexact Newton steps for a root finding problem.

\subsection{Local Contraction Property}\label{sec:analysis_contraction}
Start by formulating the first-order necessary conditions (FONC) for \eqref{eq:opt_iter}.
The FONC state that if \(\bm{y}\) is a solution to \eqref{eq:opt_iter}, then there exists a real vector \(\lambda\) of appropriate dimension such that
\begin{subequations}\label{eq:opt_iter_kkt}\begin{alignat}{1}
    2 \mathcal{Q} \bm{y} + \mathcal{G}^T(\bm{\rho}[l]) \lambda &= 0 \\
    \mathcal{G}(\bm{\rho}[l]) \bm{y} + C \hat{x}_0 &= 0
\end{alignat}\end{subequations}
hold \cite{Nocedal2006}. After introducing \(\bm{z} \coloneqq (\bm{y}, \lambda)\) as well as defining
\begin{equation}\label{eq:iter_root_function}
    F(\bm{z}) \coloneqq \bbma
        2 \mathcal{Q} \bm{y} + \mathcal{G}^T(\bm{\rho}) \lambda \\
        \mathcal{G}(\bm{\rho}) \bm{y} + C \hat{x}_0
    \ebma,
\end{equation}
observe that solving \eqref{eq:opt_iter_kkt} and
\begin{equation}\label{eq:newton_type_iteration}
    F(\bm{z}[l]) + J(\bm{z}[l]) (\bm{z} - \bm{z}[l]) = 0
\end{equation}
for \(\bm{z}\) is equivalent.
\(\bm{z}[l]\) denotes the solution of the previous iteration.
\(J(\bm{z})\) approximates \(F'(\bm{z})\), which is the Jacobian of \(F(\bm{z})\), and is given by
\begin{equation}\label{eq:iter_approximated_jacobian}
    J(\bm{z}) \coloneqq \bbma
        2 \mathcal{Q}          & \mathcal{G}^T(\bm{\rho}) \\
        \mathcal{G}(\bm{\rho}) & 0
    \ebma,
\end{equation}
while the exact Jacobian is
\begin{equation}\label{eq:iter_exact_jacobian}
    F'(\bm{z}) = \bbma
        2 \mathcal{Q} + \pd{\mathcal{G}^T(\bm{\rho})\lambda}{\bm{y}} & \mathcal{G}^T(\bm{\rho}) \\
        \pd{[\mathcal{G}(\bm{\rho})\bm{y}]}{\bm{y}}                  & 0
    \ebma.
\end{equation}
\eqref{eq:newton_type_iteration} has the structure of an approximate Newton step, thus iterating \eqref{eq:opt_iter_kkt} is equivalent to using a Newton-type scheme to solve the root finding problem \(F(\bm{z}) = 0\).
Note that \(F(\bm{z})\) is not the same function as the one that corresponds to the FONC of NLP~\eqref{eq:opt_cons}.
We will study the difference more precisely in Section~\ref{sec:analysis_optimality}.

There exists a wide body of literature for analysing Newton-type schemes which we can use to provide sufficient conditions for the iterations \eqref{eq:opt_iter} to converge.
We will only consider the case where \(\hat{x}_0\) is fixed, i.e. where \eqref{eq:opt_iter} is iterated for a single time step.
Let us make the following additional assumptions:

\begin{assumption}\label{ass:solution_regularity}
    There exists a solution point \(\bar{\bm{z}}\) of \({F(\bm{z}) = 0}\).
    Let further \(\mathcal{G}(\bar{\bm{\rho}})\) have full row rank, and \(\mathcal{Q}\) be positive definite on the nullspace of \(\mathcal{G}(\bar{\bm{\rho}})\), where \(\bar{\bm{\rho}}\) is the parameter corresponding to \(\bar{\bm{z}}\).
\end{assumption}

The assumptions made here are not very restrictive.
More precisely, if an OCP with the structure of \eqref{eq:opt_lpv_sum} is considered, \(\mathcal{G}(\bm{\rho})\) has full row rank for any \(\bm{\rho}\) and the natural choice of \(R \succ 0\) makes \(\mathcal{Q}\) positive definite on the nullspace of \(\mathcal{G}(\bm{\rho})\).
Nonetheless, they imply regularity of \eqref{eq:newton_type_iteration} around \(\bar{\bm{z}}\) as stated in the following lemma.

\begin{lemma}\label{lem:non_singularity}
    Suppose Assumptions~\ref{ass:model_differentiability} and \ref{ass:solution_regularity} hold.
    Then there exists a \(\delta > 0\) such that \(J(\bm{z})\) is \correction{non-singular} for all \(\bm{z}\) such that \(\|\bm{z} - \bar{\bm{z}}\| < \delta\).
\end{lemma}
\begin{proof}
    Assumption~\ref{ass:solution_regularity} implies that \(J(\bar{\bm{z}})\) is \correction{non-singular \cite[Lemma 16.1]{Nocedal2006}}.
    Together with Assumption~\ref{ass:model_differentiability} and continuity of the eigenvalues of a matrix, this results in \correction{the statement as given} in the lemma.
\end{proof}

Note that solving \eqref{eq:opt_iter} and \eqref{eq:opt_iter_kkt} in the vicinity of \(\bar{\bm{z}}\) is equivalent under Lemma~\ref{lem:non_singularity} since it implies that the linear independence constraint qualification holds and that all solutions to \eqref{eq:opt_iter_kkt} satisfy a second order sufficient condition and thus solve \eqref{eq:opt_iter} as well \correction{\cite[Theorem 16.2]{Nocedal2006}}.

Before presenting the main result of this section, we need to assume stronger properties than \correction{invertibility} of \(J(\bm{z})\) in a neighbourhood around \(\bar{\bm{z}}\).

\begin{assumption}\label{ass:jacobian_compatibility}
    There exist \(\omega < \infty\), \(\kappa < 1\) and an open ball \(\mathcal{B}\) around \(\bar{\bm{z}}\) such that \(J(\bm{z})\) is non-singular and
    \begin{subequations}\label{eq:local_contraction_conditions}\begin{alignat}{1}
        \|J(\bm{z})^{-1}(F'(\bm{z}) - F'(\bm{z}'))\| &\leq \omega \|\bm{z} - \bm{z}'\| \label{eq:local_contraction_conditions_lipschitz} \\
        \|J(\bm{z})^{-1}(F'(\bm{z}) - J(\bm{z}))\| &\leq \kappa \label{eq:local_contraction_conditions_compatibility}
    \end{alignat}\end{subequations}
    hold for all \(\bm{z}, \bm{z}' \in \mathcal{B}\).
\end{assumption}

In addition to the non-singularity of \(J(\bm{z})\), which is already ensured through Lemma~\ref{lem:non_singularity}, Assumption~\ref{ass:jacobian_compatibility} places two conditions on \(J(\bm{z})\) and the Jacobian \(F'(\bm{z})\): the Lipschitz condition \eqref{eq:local_contraction_conditions_lipschitz} and the compatibility condition \eqref{eq:local_contraction_conditions_compatibility}.
As noted in \cite{Diehl2016}, only condition~\eqref{eq:local_contraction_conditions_compatibility} limits for which \(F\) such a neighbourhood exists, as \eqref{eq:local_contraction_conditions_lipschitz} can, under Assumption~\ref{ass:model_differentiability}, be satisfied for any \(F\) by making \(\mathcal{B}\) and \(\omega\) sufficiently small and large, respectively.
Similar assumptions are made in \cite{Zanelli2019} for convergence analysis of SQP methods.

Now, we are ready to formulate the main result of this section, stating that \eqref{eq:newton_type_iteration}, and thus the sequence of optimization problems \eqref{eq:opt_iter}, has a local contraction property around \(\bar{\bm{z}}\).
The theorem and its proof are adapted from \cite{Diehl2016}.
Similar results that rely on slightly different assumptions can, among others, be found in \cite{Dembo1982}.

\begin{theorem}\label{thm:local_newton_contraction}
    Suppose that Assumptions~\ref{ass:model_differentiability}, \ref{ass:solution_regularity} and \ref{ass:jacobian_compatibility} hold.
    Then the sequence \(\{\bm{z}[l]\}\) generated by \eqref{eq:newton_type_iteration} converges towards \(\bar{\bm{z}}\) with contraction rate
    \begin{equation}\label{eq:local_contraction_rate}
        \left\|\bm{z}[l+1] - \bar{\bm{z}}\right\| \leq \kappa \left\|\bm{z}[l] - \bar{\bm{z}}\right\| + \frac{\omega}2 \left\|\bm{z}[l] - \bar{\bm{z}}\right\|^2
    \end{equation}
    for all \(\bm{z}[0] \in \mathcal{B}\) such that \(\|\bm{z}[0] - \bar{\bm{z}}\| < 2\frac {1 - \kappa}{\omega}\).
\end{theorem}
\begin{proof}
    Define \(\bm{z}_l \coloneqq \bm{z}[l]\), \(\Delta\bm{z}_l \coloneqq \bm{z}_l - \bar{\bm{z}}\) and \(J_l \coloneqq J(\bm{z}_l)\) and use \eqref{eq:newton_type_iteration} along with Lemma~\ref{lem:non_singularity} to see that
    \begin{equation*}
        \Delta\bm{z}_{l+1} = \Delta\bm{z}_l - J_l^{-1} F(\bm{z}_l).
    \end{equation*}
    Using the fact that \(F(\bar{\bm{z}}) = 0\), we write
    \begin{align*}
        \Delta\bm{z}_{l+1} &= \Delta\bm{z}_l - J_l^{-1} \left(F(\bm{z}_l) - F(\bar{\bm{z}})\right) \\
        &= J_l^{-1} J_l \Delta\bm{z}_l \\
        &\quad - J_l^{-1} \int_0^1 F'(\bar{\bm{z}} + t \Delta\bm{z}_l)\Delta\bm{z}_l \dif t \\
        &= J_l^{-1} \left(J_l - F'(\bm{z}_l)\right) \Delta\bm{z}_l \\
        &\quad - J_l^{-1} \int_0^1 F'(\bar{\bm{z}} + t \Delta\bm{z}_l) - F'(\bm{z}_l) \dif t \; \Delta\bm{z}_l,
    \end{align*}
    where the fundamental theorem of calculus was used for the second equality.
    Taking norms on both sides and using Assumption~\ref{ass:jacobian_compatibility}, we arrive at
    \begin{align*}
        \|\Delta\bm{z}_{l+1}\| &\leq \kappa \|\Delta\bm{z}_l\| + \omega \int_0^1 (1 - t) \dif t \; \|\Delta\bm{z}_l\|^2 \\
        &= \kappa \|\Delta\bm{z}_l\| + \frac {\omega}{2} \|\Delta\bm{z}_l\|^2.
    \end{align*}
    Convergence is then implied by \(\|\Delta\bm{z}_0\| < 2\frac{1 - \kappa}{\omega}\).
\end{proof}

Theorem~\ref{thm:local_newton_contraction} is constructive in the sense that it provides lower bounds for the convergence rate and the region of attraction of \eqref{eq:newton_type_iteration} and equivalently \eqref{eq:opt_iter}.
If only the existence of a region of attraction is of interest, other approaches will ensue less conservative results, c.f. \cite{Diehl2016}.

\subsection{Feasibility and Optimality}\label{sec:analysis_optimality}
In addition to providing conditions under which the sequence of optimization problems \eqref{eq:opt_iter} converges, we are interested in characterizing which value their solutions converge to.
The same question is answered for multiple variants of approximate SQP in \cite{Bock2007}.
The idea is to relate the stationary point of \eqref{eq:opt_iter} to a perturbed version of the NLP \eqref{eq:opt_cons}.
Before stating the resulting theorem, define
\begin{equation}\label{eq:matrix_product_rule}
    \Delta\mathcal{G}(\bm{y}) \coloneqq \pd{\big[\mathcal{G}\left(\bm{\rho}(\bm{y})\right)\bm{y}\big]}{\bm{y}}(\bm{y}) - \mathcal{G}\left(\bm{\rho}(\bm{y})\right).
\end{equation}
\(\Delta\mathcal{G}(\bm{y})\) contains the "hidden coupling" (c.f. \cite{Rugh2000}) in the quasi-LPV system caused by the dependency of \(\bm{\rho}(\bm{y})\) on \(\bm{y}\).

\begin{theorem}\label{thm:approximate_optimality}
    Suppose Assumption~\ref{ass:model_differentiability} holds, the sequence \(\{\bm{y}[l]\}\) generated by repeatedly solving QP~\eqref{eq:opt_iter} converges towards \(\bar{\bm{y}}\) for a fixed \(\hat{x}_0\), \(\bar{\lambda}\) is the corresponding multiplier at \(\bar{\bm{\rho}} \coloneqq \bm{\rho}(\bar{\bm{y}})\), and \(\bar{\bm{z}} = (\bar{\bm{y}}, \bar{\lambda})\) satisfies Assumption~\ref{ass:solution_regularity}.
    Then \((\bar{\bm{y}}, \bar{\lambda})\) satisfies the FONC of the NLP
    \begin{subequations}\label{eq:opt_perturbed_nlp}\begin{alignat}{2}
        &\!\min_{\bm{y}} & \qquad & \|\bm{y}\|_{\mathcal{Q}}^2 + e^T \bm{y} \\
        &\text{s.t.}     &        & \mathcal{G}\left(\bm{\rho}(\bm{y})\right) \bm{y} + C \hat{x}_0 = 0 \label{eq:opt_perturbed_nlp_constraints_dynamics}
    \end{alignat}\end{subequations}
    with \(e \coloneqq - \Delta\mathcal{G}^T(\bar{\bm{y}}) \bar{\lambda}\).
\end{theorem}
\begin{proof}
    Lemma~\ref{lem:non_singularity} implies constraint qualifications hold around \(\bar{\bm{y}}\), thus there exists a \(l_0 \geq 0\) such that \(\bm{y}[l]\) satisfies the FONC~\eqref{eq:opt_iter_kkt} for all \(l \geq l_0\).
    From continuity and invertibility of \(J(\bm{z})\), it follows that \(\bar{\bm{y}}\) is the unique solution of \eqref{eq:opt_iter_kkt} at \(\bar{\bm{\rho}}\), which therefore becomes
    \begin{subequations}\begin{alignat}{1}
        2 \mathcal{Q} \bar{\bm{y}} + \mathcal{G}^T(\bar{\bm{\rho}}) \bar{\lambda} &= 0 \label{eq:opt_iter_kkt_converged_lagrangian}\\
        \mathcal{G}(\bar{\bm{\rho}}) \bar{\bm{y}} + C \hat{x}_0 &= 0.
    \end{alignat}\end{subequations}
    Using \eqref{eq:matrix_product_rule}, we see that \eqref{eq:opt_iter_kkt_converged_lagrangian} is equivalent to
    \begin{equation*}
        2 \mathcal{Q} \bar{y} + e + \left[\pd{[\mathcal{G}(\bm{\rho})\bm{y}]}{\bm{y}}(\bar{\bm{y}})\right]^T \bar{\lambda} = 0,
    \end{equation*}
    such that \((\bar{\bm{y}}, \bar{\lambda})\) also satisfies the FONC of \eqref{eq:opt_perturbed_nlp}.
\end{proof}

Theorem~\ref{thm:approximate_optimality} shows that the solutions to \eqref{eq:opt_iter} do not converge to a limit that satisfies the FONC of NLP \eqref{eq:opt_cons}, but instead those of a linearly perturbed variant.
Assuming constraint qualifications hold for \eqref{eq:opt_cons}, this implies that repeatedly solving \eqref{eq:opt_iter} results in suboptimal yet feasible solutions to the NLP.

In \cite{Zanelli2019a}, it is discussed for a similar problem how the suboptimality of the solution depends on the initial condition \(\hat{x}_0\).
Under regularity assumptions on the solution of the perturbed NLP, it is shown that the suboptimality grows polynomially with the distance of \(\hat{x}_0\) from the origin, preserving stability guarantees that are formulated using the optimal solution.
A similar result for the method presented in this paper is subject to further research.

\section{Achieving Optimality}\label{sec:optimality}
Motivated by the result that iterative qLMPC converges to suboptimal solutions, we now study how to adapt the scheme in order to recover optimality.
One iterative approach which is known to converge to optimal solutions under regularity assumptions is SQP \cite{Nocedal2006}.
Applied to the NLP~\eqref{eq:opt_cons}, SQP results in a sequence of QPs given by
\begin{subequations}\label{eq:opt_sqp}\begin{alignat}{2}
    &\!\min_{\bm{p}} & \quad & \|\bm{y}[l]\|_{\mathcal{Q}}^2 + 2\bm{y}[l]^T\mathcal{Q} \bm{p} + \frac 12 \bm{p}^T \nabla_{\bm{y}}^2 \mathcal{L}[l] \bm{p} \label{eq:opt_sqp_cost}\\
    &\text{s.t.}     &       & \mathcal{G}(\bm{\rho}[l]) \bm{y}[l] + \left[\pd{[\mathcal{G}(\bm{\rho})\bm{y}]}{\bm{y}}(\bm{y}[l])\right] \bm{p} + C \hat{x}_0 = 0, \label{eq:opt_sqp_constraints_dynamics}
\end{alignat}\end{subequations}
where \(\bm{p}\) represents the solution update performed after each iteration with \(\bm{y}[l+1] = \bm{y}[l] + \bm{p}\) and \(\nabla_{\bm{y}}^2 \mathcal{L}[l]\) is the Hessian of the Lagrangian \(\mathcal{L}(\bm{y}, \lambda) \coloneqq \|\bm{y}\|_{\mathcal{Q}}^2 + \lambda^T (\mathcal{G}(\bm{\rho}(\bm{y})) \bm{y} + C \hat{x}_0)\) at the previous solution.
Since the weighted 2-norm is used as the cost function, we can apply the Gauss-Newton approximation for the Hessian, i.e. we use \(\nabla_{\bm{y}}^2 \mathcal{L}[l] \approx 2\mathcal{Q}\).
Together with a change of variables to the new optimization variable \(\bm{y}\) in the form of \(\bm{p} = \bm{y} - \bm{y}[l]\), this allows us to equivalently reformulate \eqref{eq:opt_sqp_cost} as \(\min_{\bm{y}} \|\bm{y}\|_{\mathcal{Q}}^2\).
In a second step, we define
\begin{align}\label{eq:velocity_definitions}
    \tilde{\mathcal{G}}(\bm{y}) &\coloneqq \pd{\big[\mathcal{G}\left(\bm{\rho}(\bm{y})\right)\bm{y}\big]}{\bm{y}}(\bm{y}) & \tilde{\mathcal{G}}_d(\bm{y}) &\coloneqq -\Delta\mathcal{G}(\bm{y}).
\end{align}
Using definition \eqref{eq:matrix_product_rule} on \eqref{eq:opt_sqp_constraints_dynamics} then leads to
\begin{subequations}\label{eq:opt_sqp_qlmpc}\begin{alignat}{2}
    &\!\min_{\bm{y}} & \quad & \|\bm{y}\|_{\mathcal{Q}}^2 \\
    &\text{s.t.}     &       & \tilde{\mathcal{G}}(\bm{y}[l]) \bm{y} + \tilde{\mathcal{G}}_d(\bm{y}[l]) \bm{y}[l] + C \hat{x}_0 = 0 \label{eq:opt_sqp_qlmpc_constraints_dynamics}
\end{alignat}\end{subequations}
as an equivalent optimization problem to \eqref{eq:opt_sqp}.

The structure of \eqref{eq:opt_sqp_qlmpc} is identical to that of \eqref{eq:opt_iter} with the exception of the additional \(\tilde{\mathcal{G}}_d\) term.
Notably, if an OCP in the form of \eqref{eq:opt_lpv_sum} is considered, the constraint in \eqref{eq:opt_sqp_qlmpc_constraints_dynamics} can be interpreted as the dynamics of a fictitious LPV system
\begin{equation}\label{eq:velocity_system}
    x_{k+1} = \tilde{A}(\tilde{\rho}_k) x_k + \tilde{B}(\tilde{\rho}_k) u_k + \tilde{B}_d(\tilde{\rho}_k) d_k
\end{equation}
with disturbance input \(d_k\).
This shows that one can obtain optimal solutions to the OCP using the qLMPC algorithm by implementing it for the fictitious system~\eqref{eq:velocity_system} rather than \eqref{eq:lpv_system}.
Due to the need to calculate the additional disturbance term, this leads to an increase in computational load.
The value of \(d_k\) can be calculated as a projection of \(\bm{y}[l]\), and is therefore known during optimization.

For convergence analysis of the approach described in this section we refer to comprehensive results on SQP that exist in literature, c.f. \cite{Boggs1995, Nocedal2006}.
In addition to regularity assumptions as in Section~\ref{sec:analysis}, these rely on the Hessian approximation being sufficiently good.

\section{Examples and Benchmarking}\label{sec:benchmarking}
In order to demonstrate the potential of the aforementioned methods, we will implement both on two nonlinear plants.
Additionally, two reference solvers are included in the comparison.
First, we implement the NLP~\eqref{eq:opt_cons} using \name{CasADi} \cite{Andersson2018} and solve it using \name{Ipopt} \cite{Waechter2005}.
And second, \name{acados} is included as a state-of-the-art solver in terms of computational performance.

To quantify the suboptimality introduced by each approach, we will use the notion of the relative cumulative suboptimality (RCSO) that was applied in \cite{Verschueren2019}.
The focus will be on regulation to the origin, therefore we define the distance to reference on the closed-loop solutions to the OCP as
\begin{equation}\label{eq:distance_to_reference}
    \mathrm{DR}_{(\cdot), k} \coloneqq \sum_{i=0}^k \bbma
        x_{(\cdot),i} \\ u_{(\cdot),i}
    \ebma^T \bbma
        Q & 0 \\ 0 & R
    \ebma \bbma
        x_{(\cdot),i} \\ u_{(\cdot),i}
    \ebma,
\end{equation}
where \(Q\) and \(R\) are the same matrices as used in the controller and \(k\) is the time step in the simulation.
We regard the solution generated by \name{Ipopt} as optimal and thus use it as the reference in the definition of the RCSO given by
\begin{equation}\label{eq:relative_cumulative_suboptimality}
    \mathrm{RCSO}_{(\cdot), k} \coloneqq \frac{\mathrm{DR}_{(\cdot), k} - \mathrm{DR}_{\name{Ipopt}, k}}{\mathrm{DR}_{\name{Ipopt}, k}}.
\end{equation}

For all controllers, the MPC formulation is converted into a condensed form, eliminating the state variables from the problem.
On each execution, the solvers are warm-started using the solution from the previous call.
Except for switching to Euler's method for the internal integrator, \name{acados} and \name{Ipopt} are used in their default configuration.
The simulations are run on a \SI{3.5}{\giga\hertz} Intel Core i5-4690 processor, and their source code is available at \cite{SimulationFiles}.

Since the qLMPC schemes are implemented in \name{Matlab}, the benchmarks are not fully representative of the computational performance of the approach, even though the \name{Matlab Coder} is used to generate compiled \name{MEX}-functions.
This is especially true for the second example, where the implementation of the LPV model uses multiple matrix multiplications.
In order to evaluate the intrinsic performance of the approach, an implementation in a compiled programming language should be considered in the future.

\subsection{Dynamic Unicycle}\label{sec:benchmarking_unicycle}
The first example that we consider is a dynamic unicycle.
It can be seen as a simplified model of a wheeled robot featuring non-holonomic constraints and is described by
\begin{align*}
    \dot{s}(t)    &= v(t) \cos\left(\phi(t)\right) & \dot{v}(t)      &= F(t) \\
    \dot{q}(t)    &= v(t) \sin\left(\phi(t)\right) & \dot{\omega}(t) &= \tau(t) \\
    \dot{\phi}(t) &= \omega(t).
\end{align*}
We discretize the model using Euler's method with sampling time \(T = \SI{0.1}{\second}\), which, after defining the state, input and parameter vectors as \(x_k \coloneqq (s_k, q_k, v_k, \phi_k, \omega_k)\), \(u_k \coloneqq (F_k, \tau_k)\) and \(\rho_k \coloneqq \phi_k\), leads to an LPV model with
\begin{align*}
    A(\rho_k) &= \bbma
        1 & 0 & T \cos(\phi_k) & 0 & 0 \\
        0 & 1 & T \sin(\phi_k) & 0 & 0 \\
        0 & 0 & 1              & 0 & 0 \\
        0 & 0 & 0              & 1 & T \\
        0 & 0 & 0              & 0 & 1
    \ebma &
    B &= \bbma
        0 & 0 \\
        0 & 0 \\
        T & 0 \\
        0 & 0 \\
        0 & T
    \ebma.
\end{align*}

Following the approach outlined in Section~\ref{sec:optimality}, we calculate a fictitious model in the form of \eqref{eq:velocity_system} for use in the optimality preserving qLMPC variant.
This results in
\begin{equation*}
    \tilde{A}(\tilde{\rho}_k) = \bbma
        1 & 0 & T \cos(\phi_k) & -T v_k \sin(\phi_k) & 0 \\
        0 & 1 & T \sin(\phi_k) &  T v_k \cos(\phi_k) & 0 \\
        0 & 0 & 1              & 0                   & 0 \\
        0 & 0 & 0              & 1                   & T \\
        0 & 0 & 0              & 0                   & 1
    \ebma
\end{equation*}
and \(\tilde{B} = B\), where we define \(\tilde{\rho}_k \coloneqq (v_k, \phi_k)\).
Accordingly, the disturbance input matrix is given by
\begin{equation*}
    \tilde{B}_d(\tilde{\rho}_k) = \bbma
        T v_k \sin(\phi_k) & -T v_k \cos(\phi_k) & 0 & 0 & 0
    \ebma^T
\end{equation*}
and the disturbance is \(d_k = \phi_k[l]\), i.e. the orientation that was obtained in the previous solver iteration.

The unicycle starts with initial condition \(x_0 = (1, 2, 0, \pi, 0)\) and must be steered to the origin.
We chose \(N = 20\) for the prediction horizon, \(Q = P = \diag(1, 1, 0.1, 1, 0.1)\) and \(R = \diag(1, 1)\) for the weighting matrices and the simulation is run for \SI{10}{\second}.
The qLMPC variants and \name{acados} are configured to apply a single iteration of their respective algorithms.
This corresponds to the real-time iteration scheme proposed in \cite{Diehl2005}.

\begin{figure}
    \centering
    \begin{tikzpicture}
        \begin{axis}[
            height=6cm, width=\columnwidth,
            xlabel=$s(t)$ Coordinate,
            ylabel=$q(t)$ Coordinate,
            legend pos=north west,
        ]
            \pgfplotstableread[col sep=comma]{data/dynamic_unicycle_traj.csv}\datatable

            \addplot table[x=x_qlmpc, y=y_qlmpc] {\datatable};
            \addplot table[x=x_vqlmpc, y=y_vqlmpc] {\datatable};
            \addplot table[x=x_acados, y=y_acados] {\datatable};
            \addplot table[x=x_casadi, y=y_casadi] {\datatable};

            \addplot[black, only marks, mark=x] coordinates {(0,0)};

            \legend{qLMPC, Exact qLMPC, \name{acados}, \name{Ipopt}}
        \end{axis}

        \path (current axis.south west) -- ++(0, 0) rectangle (current axis.north east) -- ++(0, 1ex);
    \end{tikzpicture}
    \caption{Closed-loop unicycle trajectories with \(N = 20\).}
    \label{fig:unicycle_trajectories}
\end{figure}
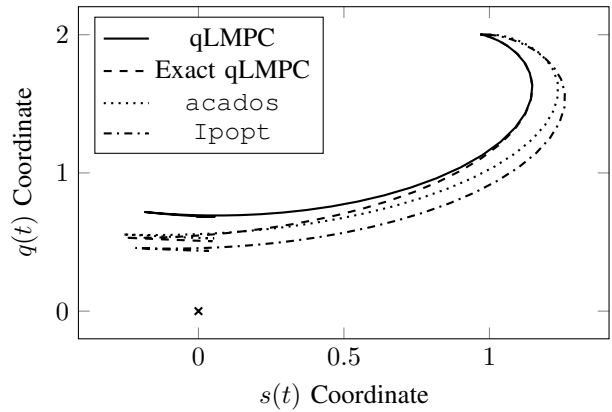

The resulting closed-loop trajectories are shown in Fig.~\ref{fig:unicycle_trajectories}.
Because of the non-holonomic constraints no controller manages to steer the unicycle to the origin.
The error is larger for standard qLMPC than for its exact variant and \name{acados} due to solution suboptimality.
This can also be seen in Table~\ref{tab:dynamic_unicycle_simulation_results}, where the RCSO is more than twice as large for standard qLMPC, albeit still having less than \SI{9}{\percent} higher cumulative cost relative to \name{Ipopt}.

The solver times in Table~\ref{tab:dynamic_unicycle_simulation_results} are obtained by running the simulation 50 times and averaging the measurements.
\begin{table}
	\centering
    \caption{Simulation Results for the Dynamic Unicycle}
	\label{tab:dynamic_unicycle_simulation_results}

	\pgfplotstableread[col sep=comma]{data/dynamic_unicycle_perf.csv} \datatable

	\pgfplotstabletypeset[
		precision=3,
		every column/.style={
			column type=r,
			fixed,
			zerofill
		},
		columns/Row/.style={
			string type,
			column name={Controller},
			column type=l,
			string replace={qlmpc}{qLMPC},
			string replace={vqlmpc}{Exact qLMPC},
			string replace={acados}{\name{acados}},
			string replace={casadi}{\name{Ipopt}}
		},
		columns/median/.style={column name=Median},
		columns/min/.style={column name=Min},
		columns/max/.style={column name=Max},
		columns/rcso/.style={column name=RCSO, sci, precision=2, clear infinite},
		every head row/.style={
			before row={
                \toprule
                & \multicolumn{3}{c}{Solver time in \si{\milli\second}} \\ \cmidrule(lr){2-4}
            },
			after row=\midrule
		},
		every last row/.style={after row=\bottomrule},
		empty cells with={---}
	]\datatable
\end{table}
It can be seen that \name{acados} and the standard qLMPC variant have identical median performance, with qLMPC falling behind if worst-case performance is considered.
The impact of using the exact qLMPC variant instead is \SI{23}{\percent} in this example.
Due to performing a single iteration in each step, all real-time methods require a constant runtime during the simulations.

\subsection{Arm-Driven Inverted Pendulum}\label{sec:benchmarking_adip}
The second model under consideration is the arm-driven inverted pendulum (ADIP) that was studied in \cite{Cisneros2019}.
Arm and pendulum are connected using a free spinning joint, and a torque~\(\tau\) can be applied to the base of the arm to control the plant.
The state of the model is given by \(x \coloneqq (\theta_1, \theta_2, \dot{\theta}_1, \dot{\theta}_2)\), the angles of the arm and the pendulum to the vertical axis and their derivatives, respectively.
As with the unicycle, we apply Euler's method for discretization.
For details on the plant and its dynamics, see \cite{Cisneros2019}.

The ADIP is initialized with \(x_0 = (\frac{\pi}{3}, 0, 0, 0)\) and should be regulated to the unstable equilibrium at the origin.
We select \(N = 40\), \(Q = P = \diag(200, 1000, 0.1, 10)\) and \(R = 2000\) and simulate the plant for \SI{2}{\second}.
For this scenario, all methods are configured to perform multiple iterations per step if required.
The qLMPC variants use the infinity norm of the FONC' residual as stopping criterion.

All controllers manage to bring the ADIP to the upright position.
As can be seen in Table~\ref{tab:adip_simulation_results}, standard qLMPC and \name{acados} achieve a RCSO between \SIrange[range-phrase ={ and }]{1}{5}{\percent}.
\begin{table}
	\centering
	\vspace{1ex}
    \caption{Simulation Results for the ADIP}
	\label{tab:adip_simulation_results}

	\pgfplotstableread[col sep=comma]{data/adip_subtask.csv} \datatable

	\pgfplotstabletypeset[
		precision=3,
		every column/.style={
			column type=r,
			fixed,
			zerofill,
		},
		columns/Row/.style={
			string type,
			column name={Controller},
			column type=l,
			string replace={qlmpc}{qLMPC},
			string replace={vqlmpc}{Exact qLMPC},
			string replace={acados}{\name{acados}},
			string replace={casadi}{\name{Ipopt}}
		},
		columns/total/.style={column name=Total},
		columns/qp/.style={column name=QP Solver, clear infinite},
		columns/prep/.style={column name=Preparation, clear infinite},
		columns/rcso/.style={column name=RCSO, sci, precision=2, clear infinite},
		every head row/.style={
			before row={
                \toprule
                & \multicolumn{3}{c}{Median solver time in \si{\milli\second}} \\ \cmidrule(lr){2-4}
            },
			after row=\midrule
		},
		every last row/.style={after row=\bottomrule},
		empty cells with={---}
	]\datatable
\end{table}
The exact qLMPC variant managed to take a better trajectory than \name{Ipopt}, which is possible due to discretization errors.

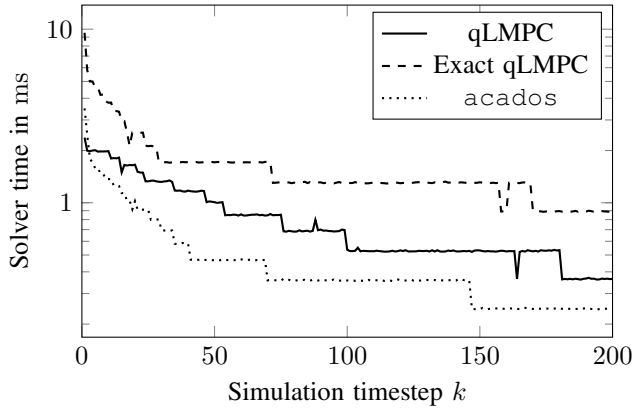
\begin{figure}
    \centering
    \begin{tikzpicture}
        \begin{semilogyaxis}[
            height=6cm, width=\columnwidth,
            xmin=0, xmax=200,
            xlabel=Simulation timestep $k$,
            ylabel=Solver time in \si{\milli\second},
            log ticks with fixed point
        ]
            \pgfplotstableread[col sep=comma]{data/adip_perf_plot.csv}\datatable

            \addplot table[x=step, y=qlmpc] {\datatable};
            \addplot table[x=step, y=vqlmpc] {\datatable};
            \addplot table[x=step, y=acados] {\datatable};

            \legend{qLMPC,  Exact qLMPC, \name{acados}}

        \end{semilogyaxis}

        \path (current axis.south west) -- ++(0, 0) rectangle (current axis.north east) -- ++(0, 1ex);
    \end{tikzpicture}
    \caption{Solver time of the three real-time methods for the ADIP with \(N = 40\), averaged over 50 runs of the simulation.}
    \label{fig:adip_solver_performance}
\end{figure}

The solver times are shown in Fig.~\ref{fig:adip_solver_performance}.
The plots go down in steps as the simulation progresses since the real-time methods perform fewer iterations as the ADIP approaches the origin.
Overall, \name{acados} outperforms qLMPC for the ADIP.
As mentioned at the beginning of this section and backed by the timings in Table~\ref{tab:adip_simulation_results}, this is at least in part due to inefficient implementation of the ADIP model.
The matrix representation of the QP is calculated in the preparation phase, which takes a large majority of the total solver time.
One approach to improve the performance could be to use \name{CasADi} to generate efficient code for the ADIP symbolically.
In this way, computationally expensive matrix inversions and multiplications could to some extent be calculated offline.

\section{Conclusions}\label{sec:conclusions}
In this paper, we studied the convergence properties of an iterative quasi-LPV MPC approach.
We showed that the scheme converges to suboptimal solutions under assumptions similar to those found in the literature for SQP methods and how to adapt the scheme to arrive at the optimal solution.
Moreover, the approach was demonstrated to be competitive with state-of-the-art solvers for nonlinear MPC in terms of computational performance in certain scenarios, while facilitating stability analysis using terminal ingredients as described in \cite{Cisneros2020}.

    \bibliographystyle{IEEETran}
    \bibliography{root}
\end{document}